\documentclass[12pt,a4paper]{article}
\usepackage{latexsym}
\usepackage{amsfonts, dsfont, eufrak, yfonts}
\usepackage{amssymb,fontenc,mathrsfs,bbm}
\newcommand{\grad}{{{\rm{grad}\,}}}

\newcommand{\Hess}{{{\rm Hess}\,}}
\newcommand{\noin}{\noindent}
\newcommand{\Hm}{{
\mathbbmss h}}

\newtheorem{theorem}{\bf Theorem}[section]

\newtheorem{remark}[theorem]{\bf Remark}

\begin{document}

\title{On compact H-Hypersurfaces of $N\times \mathbb{R}$  }
\author{G. Pacelli Bessa\thanks{Research partially supported by CNPq-Brasil} \and J. Fabio Montenegro}
\date{\today}
\maketitle
\begin{abstract}
\noin  Let ${\mathscr F}(N\times \mathbb{R})$ be the set of all
closed
 $H$-hypersurfaces $M\subset N\times \mathbb{R}$, where $N$ is a simply connected complete
 Riemannian $n$-manifold with sectional curvature
$K_{N}\leq -\kappa^{2}<0$. We  show that ${\Hm }(N\times
\mathbb{R})=\inf_{M\in {\mathscr F}(N\times \mathbb{R})}\{\vert
H_{M}\vert \}\geq (n-1)\kappa/n $.

\vspace{.2cm}


\end{abstract}
\section{Introduction}\hspace{.5cm} Let
$G(f)=\{(x,f(x))\}\subset N\times \mathbb{R}$  be a graph of a
smooth function $ f:N\to \mathbb{R}$ with constant mean curvature
$H_{G(f)}$, where $N$ is a complete Riemannian $n$-manifold. It
follows from the work of Barbosa, Kenmotsu and Oshikiri
\cite{barbosa-kenmotsu-Oshikiri} and Salavessa  \cite{salavessa2}
that
\begin{equation}\label{eqSalavessa}  \vert H_{G(f)}\vert \leq
(2/n)\cdot\sqrt{\lambda^{\ast}(N)},\end{equation} where $
\lambda^{\ast}(N)= \inf\left\{\int_{N}\vert \nabla
f\vert^{2}/\int_{N} f^{2},\,f\in
 {H_{0}^{1}(N)\setminus\{0\}}\right\}
$ is the fundamental tone of $N$ and $H_{0}^{1}(N)$ is the
completion of $C^{\infty}_{0}(N)$ with respect to the norm $
\Vert\varphi \Vert_{H_{0}^{1}}^2=\int_{N}\varphi^{2}+\int_{N}
\vert\nabla \varphi\vert^2. $\\

When $N= \mathbb{H}^{n}(-1)$ the inequality (\ref{eqSalavessa})
becomes
\begin{equation} \vert H_{G(f)}\vert \leq
(2/n)\cdot\sqrt{\lambda^{\ast}(\mathbb{H}^{n}(-1))}=(n-1)/n.
\end{equation}

On the other hand,  closed embedded $H$-hypersurface $M\subset
\mathbb{H}^{n}(-1)\times \mathbb{R}$ has mean curvature $\vert H_{M}
\vert > (n-1)/n$. This was proved by Nelli and Rosenberg
\cite{nelli} for $n=2$, independently by Salavessa
\cite{salavessa1}, \cite{salavessa2} for any $n\geq 2$, constructing
entire vertical graphs $G(f)\subset \mathbb{H}^{n}(-1)\times
\mathbb{R}$ with constant mean curvature $\vert H_{G(f)}\vert =c/n$
for each $c\in (0, n-1]$ and applying the maximum principle.
  It should be remarked that this result for $n=2$ was implicit in
Hsiang-Hsiang's paper \cite{hsiang-hsiang}. It is a consequence of
 Abresch-Rosenberg's work \cite{abresch-rosenberg}, that
any $2$-sphere  with constant mean curvature immersed in
$\mathbb{H}^{2}(-1)\times \mathbb{R}$ has mean curvature $\vert
H\vert >1/2$.

After these results, Rosenberg  \cite{rosenberg}, suggested  the
invariant ${\Hm }(N)=\inf_{M\in {\mathscr F}(N)} \{\vert
H_{M}\vert\}$, where ${\mathscr F}(N)$ is the set of all closed
 $H$-hypersurfaces $M$ immersed the complete Riemannian manifold
 $N$ and asked whether
 $\Hm
(\mathbb{H}^{n}(-1)\times \mathbb{R})\geq (n-1)/n$. The purpose of
this paper is to answer Rosenberg's question affirmatively. In fact,
we prove the following more general result.

\begin{theorem}\label{thm1}Let $N$ be a   complete
$n$-dimensional Riemannian manifold with a pole  and radial
sectional curvature bounded above\footnote{Meaning that  the
sectional curvatures along the geodesics emanating from the pole $q$
are bounded above.} $K_{N}\leq -\kappa^{2}< 0$, $($$\kappa>0$$)$.
Let $M\subset N\times \mathbb{R}$ be a closed immersed
H-hypersurface with  mean curvature $H$. Then
\begin{equation}\label{eqBF0} \vert H_{M}\vert \geq
(2/n)\sqrt{\lambda^{\ast}(\mathbb{H}^{n}(-\kappa^{2}))}=(n-1)\cdot
\kappa/n.
\end{equation}In particular \begin{equation}\Hm ( N\times
\mathbb{R})\geq (n-1)\cdot \kappa/n.
  \end{equation}
\end{theorem}

\begin{remark}The  geodesic spheres $\partial
B_{\mathbb{R}^{n}}(r)\subset\mathbb{R}^{n}$ of radius $r$ has
constant mean curvature $\vert H_{\partial B(r)}\vert =1/r$.
Therefore $\Hm (\mathbb{R}^{n})=0$. The totally
 geodesics $(n-1)$-sphere $S^{n-1}(1)\subset S^{n}(1)$ shows that $\Hm (\mathbb{S}^{n}(1))=0$.
   Jorge and Xavier
\cite{jorge-xavier} showed $($in particular$)$ that  $\Hm(
\mathbb{H}^{n}(-1))\geq 1$. The geodesic sphere $\partial
B_{\mathbb{H}^{n}(-1)}(r)\subset\mathbb{H}^{n}(-1)$  has constant
mean curvature $\vert H_{\partial B(r)}\vert=\coth(r)$. Thus we have
that  $\Hm (\mathbb{H}^{n}(-1))= 1$.
\end{remark}

 Let
 $p_{1}: N\times \mathbb{R}\to N$ the projection on the first factor
 and for a
 given closed hypersurface $M\subset N\times \mathbb{R}$ let
 $R_{p_{1} (M)}$ be the extrinsic radius of the set $p_{1} (M)\subset N$.
 Our second result shows that if there is a sequence of closed immersed H-hypersurfaces $M_{i}\subset N\times
\mathbb{R}$ with constant mean curvature $\vert H_{M_{i}}\vert
\to(n-1)\cdot\kappa/n$ then  $R_{p_{1} (M_{i})}\to \infty$. We prove
the following theorem.
\begin{theorem}\label{thm2}Let $N$ be a   complete
$n$-dimensional Riemannian manifold with a pole  and radial
sectional curvature bounded above $K_{N}\leq -\kappa^{2}< 0$,
$($$\kappa>0$$)$. Let $M\subset N\times \mathbb{R}$ be a closed
immersed hypersurface with constant mean curvature $H_{M}$. Then the
extrinsic radius
\begin{equation}R_{p_{1}(M)}\geq \coth^{-1}(\frac{n }{(n-1)\cdot
\kappa}\cdot \vert H_{M}\vert )\end{equation}In particular if
$M_{i}\subset N\times \mathbb{R}$ is a sequence of  closed immersed
hypersurfaces with constant mean curvatures $\vert H_{M_{i}}\vert
\to (n-1)\cdot\kappa/n$ then the extrinsic radius
$R_{p_{1}(M_{i})}\to \infty$.
\end{theorem}In \cite{frankel}, Frankel proved  that any two minimal
hypersurfaces, one closed and the  other properly immersed in a
complete Riemannian manifold with positive Ricci curvature  must
intersect. This result has a version for $N \times \mathbb{R}$ where
$N$ has positive sectional curvature.
\begin{theorem}\label{thm3}Let $N$ be  complete Riemannian manifold with
positive sectional curvature. Let $M_{1}$ and $M_{2}$ be two minimal
immersed hypersurfaces immersed in $N\times \mathbb{R}$, where
$M_{1}$ is closed and $M_{2}$ is proper. Then $($up to a vertical
translation$)$ they intersect, i.e. $p_{1}(M_{1})\cap
p_{1}(M_{2})\neq \emptyset$.
\end{theorem}

\section{Preliminaries}
Let $\varphi : M \hookrightarrow N$ be an isometric immersion, where
$M$ and $N$ are complete Riemannian
 manifolds. Consider a smooth function $g:N \rightarrow \mathbb{R}$ and the composition $f=g\,\circ\,
 \varphi :M \rightarrow \mathbb{R}$. Identifying $X$ with $d\varphi (X)$ we have at $q\in M$ and  for every
 $X\in T_{q}M$ that
\begin{eqnarray} \langle \grad\,f\, , \, X\rangle \,=\,d f(X)=d g(X)\,=\, \langle \grad\,g\, ,\, X\rangle,\nonumber
\end{eqnarray}
therefore
\begin{equation} \grad\, g\,=\, \grad\, f\, +\,(\grad\,g)^{\perp},\label{eqBF1}
\end{equation}
where $(\grad\,g)^{\perp}$ is perpendicular to $T_{q}M$. Let $\nabla
$ and $\overline{\nabla}$ be the Riemannian connections on $M$ and
$N$ respectively,  $\alpha (q) (X,Y) $ and $\Hess\,f(q)\,(X,Y)$  be
respectively the second fundamental form of the immersion $\varphi $
and
 the Hessian of $f$ at $q\in M$,  $X,Y \in T_{p}M$. Using the Gauss equation we have that
 \begin{equation}\Hess\,f (q)  \,(X,Y)= \Hess\,g (\varphi (q))\,(X,Y) +
 \langle \grad\,g\,,\,\alpha (X,Y)\rangle_{\varphi (q)}.
\label{eqBF2}
\end{equation}
 Taking the trace in (\ref{eqBF2}), with respect to an orthonormal basis $\{ e_{1},\ldots e_{m}\}$
 for $T_{q}M$, we have  that
\begin{eqnarray}
\Delta \,f (p)                & = & \sum_{i=1}^{m}\Hess\,f (q)  \,(e_{i},e_{i})\nonumber \\
                             & = & \sum_{i=1}^{m}\Hess\,g (\varphi (q))\,(e_{i},e_{i}) + \langle \grad\,g\,,\,
                             \sum_{i=1}^{m}\alpha (e_{i},e_{i})\rangle.\label{eqBF3}
\end{eqnarray}
 We should mention that the formulas (\ref{eqBF2}) and (\ref{eqBF3}) are  well known in the literature,
 see \cite{bessa-montenegro1}, \cite{bessa-jorge-lima-montenegro},
  \cite{choe-gulliver}, \cite{cheung-leung}, \cite{jorge-koutrofiotis}. Another important tool  is the Hessian Comparison Theorem.

\begin{theorem}[Hessian Comparison Thm.] Let $M$ be    a  complete Riemannian $n$-manifold  and
$x_{0},x_{1} \in M $.   Let
 $\gamma:[0,\,\rho (x_{1}) ]\rightarrow M$ be a minimizing geodesic joining $x_{0}$ and $x_{1}$ where $\rho (x)$
 is the  distance function $dist_{M}(x_{0}, x) $. Let $K_{\gamma}$ be the sectional curvatures of $M$ along $\gamma$
   and  let $\mu(\rho )$
  be this function defined below.

\begin{equation} \mu(\rho )=\left\{ \begin{array}{lcll}
 & k \cdot\coth (k \cdot\rho (x)), & if  & \sup K_{\gamma}=-k^{2} \\
 &                       &     &   \\
 &\displaystyle \frac{1}{\rho (x)},   & if  &  \sup K_{\gamma}=0  \\
 &                       &     & \\
 & k\cdot\cot (k\cdot \rho (x)),  & if  &  \sup K_{\gamma} =k^{2}\; and \; \rho < \pi/2k.
\end{array}\right.\label{eqBF5}
\end{equation}

\noin Then the Hessian of  $\rho$ and $\rho^{2}$  satisfies

\begin{equation}\begin{array}{lllllll}

Hess\,\rho(x)(X,X)&\geq & \mu(\rho(x))\cdot\Vert X\Vert^{2},
                                      &      & Hess\,\rho^{2}(x)(X,X)&\geq &
                                      2\rho(x)\cdot\mu(\rho(x))\cdot\Vert X\Vert^{2}\\
                                      &&&&&&\\
Hess\,\rho(x) (\gamma ',\gamma ') & =   & 0, & & Hess\,\rho^{2}(x)
(\gamma ',\gamma ')
                  & =   & 2.
\end{array}\label{eqBF6}
\end{equation}
Where $X$ is any vector in $T_{x}M$ perpendicular to $\gamma'(\rho
(x))$. \label{thmHess}

\end{theorem}

\section{Proof of the Results}
  The Theorems (\ref{thm1}) and (\ref{thm2}) are
consequences of this  following result.
\begin{theorem}\label{thmGeral} Let  $\varphi
:M\hookrightarrow N \times \mathbb{R}$ be an $n$-dimensional closed
immersed hypersurface $M$ with bounded mean curvature $\vert {H}_{M}
\vert \leq H_{0}$, where $N$ is  complete Riemannian manifold  with
a pole $x_{o}$ and
 radial sectional curvatures
bounded  $K_{N}\leq c$ and $\rho_{N}:N\to\mathbb{R}$ be the distance
function to  $x_{o}$. Suppose  that $p_{1}(M)\subset
B_{N}(\pi/2\kappa)$ if $c=\kappa^{2}$. Then
\begin{equation}n\cdot \vert H_{0}\vert \geq (n-1)\cdot \mu (R_{p_{1}(M)})
\end{equation}

\end{theorem}

\subsection{Proof of Theorem \ref{thmGeral}}

To prove Theorem \ref{thmGeral} we proceed as follows: let
$g:N\times \mathbb{R} \to \mathbb{R}$ be given by
$g(x,t)=\rho^{2}_{N}(x)$ and let $f:M\to \mathbb{R}$ be defined by
$f=g\circ \varphi$. The function $g$ is smooth on $N\times
\mathbb{R}$, and $f\geq 0$ is smooth on $M$. By Green's theorem we
have that $\smallint_{M}[f\triangle f + \vert \grad f\vert^{2}]=0$.
This implies that there exists a subset $\emptyset \neq S\subset M$
such that for any point $q\in S$ we have that $\triangle f (q) <0$.
By (\ref{eqBF3}) we have that
\begin{equation}\label{eqBF6} 0>\triangle f(q)=
\sum_{i=1}^{n}\Hess\,g (\varphi (q))\,(e_{i},e_{i}) +
                              \langle \grad\,g\,,\,
                             \sum_{i=1}^{n}\alpha
                             (e_{i},e_{i})\rangle.\end{equation}
 We choose   an orthonormal basis  $\{
e_{1},\ldots e_{n}\}$
 for $T_{q}M$ in the following way.
  Let start with an orthonormal basis (from polar coordinates) for $T_{p(q)}N$,
  $\{\grad \rho_{N},\partial/\partial \theta_{2},\ldots,
  \partial/\partial \theta_{n}\}$. We can choose be a orthonormal  basis for $T_{q}M$ as
  follows
  $e_{1}=\langle e_{1}, \partial /\partial t\rangle \partial /\partial t
 + \langle e_{1}, \grad \rho_{N}\rangle\grad \rho_{N}$ and $e_{j}=\partial/\partial \theta_{j}$, $j=2,...,n$
 (up to an re-ordination). Computing
 $\Hess g (e_{i},e_{i})$ we obtain that
 \begin{equation}\Hess g (e_{i},e_{i})=\left\{\begin{array}{lll}2\langle e_{1}, \grad \rho_{N}\rangle^{2}&
 if & i=1\\
 && \\
 2\rho_{N}\Hess\rho_{N}(e_{i},e_{i}) &if & i\geq 2
 \end{array}\right.
 \end{equation}Therefore at $x_{o}\neq q\in S$ we have that
 \begin{eqnarray}\label{eqBF7}0&>&2\langle e_{1}, \grad \rho_{N}\rangle^{2}
 +2\cdot(n-1)\cdot\rho_{N}\cdot\Hess\rho_{N}(e_{i},e_{i})+ 2\cdot  \rho_{N} \langle \rho_{N}\,,\,
                             \stackrel{\rightarrow}{ H_{M}}\rangle\nonumber \\
                             &&\nonumber \\
                             &\geq&  2\langle e_{1}, \grad \rho_{N}\rangle^{2}
                             +2\cdot(n-1)\cdot\rho_{N}\cdot\Hess\rho_{N}(e_{i},e_{i})-2\cdot n \cdot \rho_{N}\cdot\vert H_{0}\vert.
                             \\
                             && \nonumber \\
                             & \geq &
                             2\cdot(n-1)\cdot\rho_{N}\cdot\mu(\rho_{N})-2\cdot n\cdot \rho_{N}\cdot\vert H_{0}\vert\nonumber
 \end{eqnarray}From (\ref{eqBF7}) we obtain $$n\cdot \vert H_{0}\vert\geq
 (n-1)\cdot \mu (\rho_{N})\geq(n-1)\cdot \mu (R_{p(M)}).$$Observe
 that for us, $\stackrel{\rightarrow}{H_{M}}=\sum_{i=1}^{n}\alpha
                             (e_{i},e_{i})$ so that
                             $\vert\stackrel{\rightarrow}{H_{M}}\vert =
                             n\cdot \vert H_{M}\vert$.

If $\varphi:M\hookrightarrow N\times \mathbb{R}$ is a closed
H-hypersurface with mean curvature $\vert H_{M}\vert$, where $N$ has
sectional curvature $K_{N}\leq -\kappa^{2}<0$ then
\begin{equation}\label{eqBF8}n\cdot \vert H_{M}\vert \geq
(n-1)\cdot\kappa \cdot\coth(\kappa\cdot R_{p(M)})\geq (n-1)\cdot
\kappa\end{equation} This proves (\ref{eqBF0}). Now from
(\ref{eqBF8}) we can conclude that
$$R_{p(M)}\geq (1/\kappa)\cdot
\coth^{-1}(n\vert H_{M}\vert /(n-1)\kappa).$$

\subsection{Proof of Theorem \ref{thm3}}
The proof of Theorem (\ref{thm3})  is just an observation on the
proof of Frankel's Generalized Hadamard Theorem \cite{frankel}.  We
will present  his proof and make the due observation. Let $M_{1}$
and $M_{2}$ be minimal hypersurfaces of $N\times \mathbb{R}$, where
$N$ has positive Ricci curvature, $M_{1}$ is closed and $M_{2}$ is
proper. Suppose that $p_{1}(M_{1})\cap p_{1}(M_{2})=\emptyset $,
recalling that $p_{1}:N\times \mathbb{R}\to N$ is the projection on
the first factor. Let $\gamma$ be a geodesic joining $a\in M_{1}$
and $b\in M_{2}$ of positive length realizing the distance $l$
between $M_{1}$ and $M_{2}$. This geodesic hits $M_{1}$ and $M_{2}$
perpendicularly at $a$ and $b$ respectively. Let $X(0)\in
T_{a}M_{1}$ be a unit vector and $X(t)$ its parallel transport along
$\gamma$. This vector gives rise to a variation of $\gamma$ keeping
the end-points on  $M_{1}$ and $M_{2}$. The second variation formula
of the arc-length gives $L_{X}''(0)=\alpha_{2}
(X(l),X(l))-\alpha_{1}(X(0),X(0)) - \int_{0}^{l}K(X\wedge T)dt$,
where $\alpha_{1}$,$\alpha_{2}$ are the second fundamental forms of
$M_{1}$ and $M_{2}$ at  $a$ and $b$ evaluated at the vectors $X(0)$
and $X(l)$. Taking an orthonormal basis $\{X^{1}, ..., X^{n}\}$ of
$T_{a}M_{1}$ and adding up the second variation formulas we obtain
that $\sum_{i=1}^{n}L_{X^{i}}''(0)=-\int_{0}^{L}Ric_{N\times
\mathbb{R}} (\gamma'(t))dt$. Clearly
$\gamma'(t)=\gamma'_{N}(t)+\gamma'_{\mathbb{R}}(t)$ has horizontal
component $\gamma'_{N}(t)\neq 0$ for every $t\in I$ in a positive
measure subset $I \subset [0, l]$. Computing $Ric_{N\times
\mathbb{R}}(\gamma'(t))=\vert \gamma'_{\mathbb{R}}\vert
\sum_{i=1}^{n}K_{N}[(X^{i}_{N}/\vert X^{i}_{N}\vert)\wedge
(\gamma'_{\mathbb{R}}/\vert \gamma'_{\mathbb{R}}\vert)]\vert
X^{i}_{N}\vert^{2}>0$ if $K_{N}>0$. Therefore
$\sum_{i=1}^{n}L_{X^{i}}''(0)<0$ contradicting that fact that
$\gamma$ was of minimal length. This proves Theorem (\ref{thm3}).\\

\noin {\bf Acknowledgments:} We would like to thank Professor H.
Rosenberg and our friends L.
Jorge and J. de Lira for many helpful discussions on this paper.\\

\noin {\it Address of the authors:}\\
\noin Departamento de Matematica\\
\noin Campus do Pici, Bloco 914\\
\noin Universidade Federal do Cear\'{a}-UFC\\
\noin 60455-760 Fortaleza-Cear\'{a}\\
\noin Brazil\\

\noin  bessa@math.ufc.br \& fabio@mat.uf.br

\end{document}